\documentclass[11pt,a4paper,lualatex]{amsart}

\usepackage[marginparwidth=0pt,margin=20truemm]{geometry} % margin

\usepackage{mypackage} % My packages.
\usepackage{mycommand} % My commands.
\usetikzlibrary{knots}
\usetikzlibrary{positioning}

\DeclareMathOperator{\WRT}{WRT}
\DeclareMathOperator{\CS}{CS}
\DeclareMathOperator{\rmN}{N} %CGP invariant

\renewcommand{\Vert}{\mathrm{Vert}}

\usepackage[pdfencoding=auto,hypertexnames=false]{hyperref}
\hypersetup{colorlinks=false}

\numberwithin{equation}{section} % Setting of equation numbers 
\usepackage{mytheoremeng} % My theorem environments (English) with cleveref package

\begin{document}
% --------------------------------------------------------------------------

\title[$ L $-function invariants and generalized Bernoulli polynomials]{$ L $-function invariants for 3-manifolds and relations between generalized Bernoulli polynomials}

\author[Y. Murakami]{Yuya Murakami}
\address{Faculty of Mathematics, Kyushu University, 744, Motooka, Nishi-ku, Fukuoka, 819-0395, Japan}
\email{murakami.yuya.896@m.kyushu-u.ac.jp}

\date{\today}

\maketitle

% --------------------------------------------------------------------------

\begin{abstract}
	We introduce $ L $-functions attached to negative definite plumbed manifolds as the Mellin transforms of homological blocks.
	We prove that they are entire functions and their values at $ s=0 $ are equal to the Witten--Reshetikhin--Turaev invariants by using asymptotic techniques developed by the author in the previous papers.
	We also prove that linear relations between special values at negative integers of some $ L $-functions, which are common generalizations of Hurwitz zeta functions, Barnes zeta functions and Epstein zeta functions.
\end{abstract}

% --------------------------------------------------------------------------

\tableofcontents

% --------------------------------------------------------------------------

\section{Introduction} \label{sec:intro}

% --------------------------------------------------------------------------

Gukov--Pei--Putrov--Vafa~\cite{GPPV} introduced $ q $-series invariants of $ 3 $-manifolds which are called the Gukov--Pei--Putrov--Vafa (GPPV) invariants, $ \widehat{Z} $-invariants or homological blocks.
One of the most important properties of GPPV invariants is that their radial limits coincide with the Witten--Reshetikhin--Turaev (WRT) invariants~\cite{Reshetikhin-Turaev,Witten}, which are quantum invariants of $ 3 $-manifolds.
Such a property is conjectured by Gukov--Pei--Putrov--Vafa~\cite{GPPV} and proved by Murakami~\cite{M:plumbed,M:GPPV} by developing two methods; comparison of asymptotic expansions and pruning of plumbed graphs.
The first one is based on the technique to compute asymptotic expansions by use of the Euler--Maclaurin summation formula \cite{BKM,BMM:high_depth,M:plumbed,Zagier:asymptotic}.
On the other hand, we can compute asymptotic expansions by using Mellin transform \cite{Lawrence-Zagier,Zagier:asymptotic}.
In this article, we combine these two techniques to compute asymptotic expansions.
Consequently, we prove that WRT invariants can be written as the special values of $ L $-functions defined as the Mellin transforms of GPPV invariants and such $ L $-functions are entire functions.
Moreover, by applying this method to a general setting independent of topology, we also prove linear relations of special values at negative integers of some $ L $-functions which are common generalizations of Hurwitz zeta functions, Barnes zeta functions and Epstein zeta functions.
Here, we remark that such special values are generalizations of Bernoulli polynomials.

% --------------------------------------------------------------------------

\subsection{$ L $-function invariants} \label{subsec:intro_L-func_inv}

% --------------------------------------------------------------------------

Let $ M $ be a rational homology sphere, that is, a closed oriented 3-manifold with $ H_1(M, \Q) = 0 $.
For a positive integer $ k $, denote by $ Z_k (M) \in \bbC $ the $ \SU(2) $ Witten--Reshetikhin--Turaev (WRT) invariants at the level $ k $ of $ M $.%
\footnote{In this article, we normalize the WRT invariants as $ Z_k(S^3) = 1 $.}
In the study of the WRT invariants, there are two important open problems; construction of homology theory categorifying the WRT invariants and the asymptotic expansion conjecture~\cite{Andersen-Himpel-Jorgensen-Martens-McLellan} which conjectured expressions of asymptic expansions of WRT invariants $ Z_k (M) $ as $ k \to \infty $ in terms of the Chern--Simons invariants.
To prove these problems, the WRT invariants would need to be expressed as radial limits of some $ q $-series.
Such a strategy is established firstly by Lawrence--Zagier~\cite{Lawrence-Zagier} to prove the asymptotic expansion conjecture for the Poincar\'{e} homology sphere and developed by Hikami~\cite{Hikami:Bries,Hikami:Lattice,Hikami:Lattice2,Hikami:Seifert}, Fuji--Iwaki--Murakami--Terashima~\cite{FIMT}, Andersen--Misteg\aa{}rd~\cite{Andersen-Mistegard} and Matsusaka--Terasima~\cite{Matsusaka-Terashima} for Seifert homology spheres.
Gukov--Putrov--Vafa~\cite{GPV} and Gukov--Pei--Putrov--Vafa~\cite{GPPV} systematically interpreted such $ q $-series.
They provided a physical construction of homology theory categorifying the WRT invariants and defined $ q $-series invariants as their Euler characteristics.
Moreover, Gukov--Pei--Putrov--Vafa~\cite{GPPV} defined $ q $-series invariants mathematically for negative definite plumbed manifolds, which form a subclass of rational homology spheres.
Their $ q $-series invariants are called the Gukov--Pei--Putrov--Vafa (GPPV) invariants, $ \widehat{Z} $-invariants or homological blocks.

For a general rational homology sphere, the existence of $ q $-series whose radial limits are the WRT invariants is called the radial limit conjecture and is formulated as follows. 

\begin{conj}[{Radial limit conjecture, Hikami~\cite[Equation (1.4)]{Hikami:radial_limit}, Gukov--Pei--Putrov--Vafa~\cite[Conjecture 2.1]{GPPV}, Gukov--Manolescu~\cite[Conjecture 3.1]{Gukov-Manolescu}}]
	\label{conj:radial_limit_su(2)} 
	Let $ \Spin^c(M) $ be the set of $ \Spin^c $ structure of $ M $ acted by $ \{ \pm 1 \} $ via the conjugation.
	Then, for each $ b \in \Spin^c(M) / \{ \pm 1 \} $, there exist topological invariants
	\[
	\Delta_b \in \Q, \quad
	c \in \Z_{\ge 0}, \quad
	\widehat{Z}_b (q; M) \in 2^{-c} q^{-\Delta_b} \Z[[q]]
	\]
	such that $ \widehat{Z}_b (q; M) $ converges on $ \abs{q} < 1 $ and for infinitely many $ k \in \Z_{>0} $ it holds that
	\[
	Z_k (M) 
	=
	\lim_{q \to \zeta_k}
	\frac{1}{\zeta_{2k} - \zeta_{2k}^{-1}}
	\sum_{a, b \in \Spin^c(M) / \{ \pm 1 \}} \bm{e} \left( k \lk(a, a) \right) S_{a b} \widehat{Z}_b (q; M),
	\]
	where
	\begin{itemize}
		\item  $ \zeta_k \coloneqq e^{2\pi\iu/k} $;
		\item $ \lk \colon H_1(M,\Z) \otimes H_1(M,\Z) \to \Q/\Z $ is the linking form;
		\item $ \Stab_{\{ \pm 1 \}} (a) $ is the stable subgroup of $ a \in \Spin^c(M) $;
		\item $ \bm{e}(z) \coloneqq e^{2\pi\iu z} $ for a complex number $ z $,
		\item 	\[
		S_{a b} \coloneqq
		\frac{e^{2\pi\iu k \lk(a, b)}+ e^{-2\pi\iu k \lk(a, b)}}{\abs{\Stab_{\{ \pm 1 \}} (a)} \abs{\Stab_{\{ \pm 1 \}} (b)} \sqrt{\abs{H_1(M, \Z)}}}.
		\]
	\end{itemize}
\end{conj}

For negative definite plumbed manifolds $ M $, Gukov--Pei--Putrov--Vafa~\cite[Subsection 3.4]{GPPV} constructed the $ q $-series invariants $ \widehat{Z}_b (q; M) $ in \cref{conj:radial_limit_su(2)} called GPPV invariants as we already mentioned.
Murakami~\cite[Theorem 1.2]{M:GPPV} proved \cref{conj:radial_limit_su(2)} on the basis of studies~\cite{BMM:high_depth,FIMT,Lawrence-Zagier,MM,M:plumbed}.
For plumbed homology spheres determined by indefinite H-graphs, which is the simplest case of not negative definite, Murakami~\cite{M:indefinite} constructed $ q $-series as indefinite false theta functions and proved \cref{conj:radial_limit_su(2)}.

Under \cref{conj:radial_limit_su(2)}, we introduce the $ L $-function invariant for 3-manifold $ M $ indexed by $ b \in \Spin^c(M) / \{ \pm 1 \} $ as
\[
L_b (s; M, k)
\coloneqq
\frac{1}{\Gamma(s)} 
\int_{0}^{\infty} \widehat{Z}_b \left( \zeta_k e^{-t}; M \right) t^{s-1} dt.
\]
For negative definite plumbed manifolds or plumbed homology spheres determined by indefinite H-graphs, these $ L $-functions exist without assuming \cref{conj:radial_limit_su(2)} by the construction of GPPV invariants~\cite{GPPV} or \cite{M:indefinite} respectively.
Since Gukov--Manolescu~\cite[Proposition 4.6]{Gukov-Manolescu} proved that Gukov--Pei--Putrov--Vafa invariants $ \widehat{Z}_b (q; M) $ are topological invariants of negative definite plumbed manifolds, $ L_b (s; M, k) $ is so.
We will express $ L_b (s; M, k) $ explicitly for negative definite plumbed manifolds in \cref{subsec:proof_WRT}.

The following statement is our first main result.

\begin{thm} \label{thm:L-func_WRT}
	Let $ M $ be a negative definite plumbed manifold or plumbed homology spheres determined by indefinite H-graphs and $ b \in \Spin^c(M) / \{ \pm 1 \} $. 
	Then, the $ L $-function $ L_b (s; M, k) $ extends meromorphically on $ \bbC $.
	Moreover, the $ L $-function
	\[
	\sum_{a, b \in \Spin^c(M) / \{ \pm 1 \}} \bm{e} \left( k \lk(a, a) \right) S_{a b} L_b (s; M, k)	
	\]
	is entire and its special value at $ s=0 $ is equal to $ (\zeta_{2k} - \zeta_{2k}^{-1}) Z_k(M) $.
	These hold for general rational homology spheres $ M $ under \cref{conj:radial_limit_su(2)}.
\end{thm}

It is expected that $ \widehat{Z}_b(q; M) $ in \cref{conj:GPPV} also express non-semisimple quantum invariants called the Costantino--Geer--Patureau-Mirand invariant.
Let $ M $ be a rational homology sphere, $\omega \in H^1(M, \Q/2\Z) \smallsetminus H^1(M, \Z/2\Z)$ and $ r \in \Z_{\ge 2} \smallsetminus 4\Z$.
Denote the Costantino--Geer--Patureau-Mirand (CGP) invariant by $ \rmN_r(M,\omega) \in \bbC$, which was defined in \cite{CGP} by using the non-semisimple ribbon tensor category associated with the so-called ``unrolled'' quantum group $U_{\zeta_{2r}}^H(\mathfrak{sl}(2,\bbC)) $ constructed in \cite{CGP}.

Costantino--Gukov--Putrov established the radial limit conjecture for the CGP invariant as follows.

\begin{conj}[{\cite[Conjecture 1.2]{Costantino-Gukov-Putrov}}] \label{conj:Costantino-Gukov-Putrov}
	Under the above settings and \cref{conj:GPPV}, it holds that
	\[
	\rmN_r (M,\omega)
	=
	\lim_{q \to \zeta_r}
	\sum_{b \in \Spin^c(M)} c^{\mathrm{CGP}}_{\omega, b} \widehat{Z}_b (q; M)		
	\]
	with
	\begin{align} \label{coeffsummary}
		c^{\mathrm{CGP}}_{\omega,\sigma(b,s)}
		= \,
		&\frac{\mathcal{T}(M, [\omega])}{\abs{H_1(M;\Z)}}
		\\
		&\cdot
		\begin{dcases}
			-\bm{e} \left( -\frac{\mu(M,s)}{4} \right)
			\sum_{a, f \in H_1(M, \Z)}
			\bm{e} \left( -\frac{r-1}{4}\lk(a,a) +\lk(a,f-b) - \frac{1}{2}\omega(a) +\lk(f,f) \right)
			&\text{ if } r \equiv 1\bmod 4,
			\\
			\sqrt{|H_1(M;\Z)|} \sum_{a \in H_1(M, \Z)}
			\bm{e} \left( -\frac{q_s(a)}{4} - \lk(a,b) - \frac{\omega(a)}{2} \right),
			&\text{ if } r \equiv 2\bmod 4,
			\\
			\bm{e} \left( \frac{\mu(M,s)}{4} \right)
			\sum_{a, f \in H_1(M, \Z)}
			\bm{e} \left( -\frac{r+1}{4}\lk(a,a) -\lk(a,f+b)-\frac{1}{2}\omega(a) -\lk(f,f) \right)
			&\text{ if } r \equiv 3\bmod 4,
		\end{dcases}
	\end{align}
	where
	\begin{itemize}
		\item $\mathcal{T}(M,[\omega])$ is the Reidemeister torsion;
		\item $ \sigma \colon \ H_1(M, \Z) \times \Spin(M) \to \Spin^c(M) $ is the canonical map:
		\item $\mu(M,s)$ is the Rokhlin invariant of $M$ for spin structure $ s \in \Spin(M) $.
	\end{itemize}
\end{conj}

For negative definite plumbed manifolds $ M $, Misteg\aa{}rd--Murakami~\cite[Theorem 1.2]{M-Mistegard} proved \cref{conj:Costantino-Gukov-Putrov}.
Consequently, we obtain the following result.

\begin{thm} \label{thm:L-func_CGP}
	Let $ M $ be a negative definite plumbed manifold. 
	Then, the $ L $-function
	\[
	\sum_{b \in \Spin^c(M)} c^{\mathrm{CGP}}_{\omega, b} L_b (s; M, r)	
	\]
	is entire and its special value at $ s=0 $ is equal to the CGP invariant $ \rmN_r (M,\omega) $.
	These hold for general rational homology spheres $ M $ under \cref{conj:Costantino-Gukov-Putrov}.
\end{thm}

% --------------------------------------------------------------------------

\subsection{Functional equation} \label{subsec:intro_functional_equation}

% --------------------------------------------------------------------------

The functional equation of the Riemann zeta function can be derived from the modular transformation law of the theta function.
A similar approach can be applied to our $ L $-function.
For Seifert homology spheres, the modular transform of GPPV invariants is essentially obtained by Andersen--Misteg\aa{}rd~\cite[p.751]{Andersen-Mistegard}.

For simplicity, we only consider the $ L $-function
\[
L_0 (s; M)
\coloneqq
\frac{1}{\Gamma(s)} 
\int_{0}^{\infty} \widehat{Z}_0 \left( e^{-t}; M \right) t^{s-1} dt.
\]
Using the above modular transformation formula, we can drive the functional equation of this $ L $-function as follows.

\begin{thm} \label{thm:functional_equation}
	For a Seifert homology sphere $ M $ with $ n $ exceptional fibers and $ s \in \bbC $ with $ \Re(s) > -1 $, it holds that
	\begin{align}
		(-1)^n \pi \Gamma(s) L_0 (s; M)
		&=
		\sum_{d=0}^{n-2}
		\left( 2\pi\iu \right)^{3/2-d}
		\left( \frac{\pi}{\sqrt{\Delta P}} \right)^{s+d-1/2}
		\sum_{m=1}^{\infty}
		g_d (m) m^{s+d-1/2} K_{s+d-1/2} \left( 2\pi \sqrt{\frac{\Delta}{P}} m \right)
		\\
		&\phant
		+ 4\sqrt{\pi} \Delta^{-s/2+1/4}
		\int_{0}^{e^{- \varepsilon \iu} \infty}
		K_{s-1/2} (2\sqrt{\Delta}x)
		G \left( \exp \left( -\frac{x\iu}{\sqrt{P}} \right) \right)
		x^{s - 1/2} dx,
	\end{align}
	where $ P, \Delta \in \Z_{>0} $, $ g_0, \dots, g_{n-2} \colon \Z_{>0} \to \bbC $ and $ G(z) \in \Q(z) $ are determined by $ M $ as in \cref{subsec:proof_functional_equation}, $ K_s(z) $ is the $ K $-Bessel function and $ \varepsilon > 0 $ is a sufficiently small number.
\end{thm}

% --------------------------------------------------------------------------

\subsection{$ L $-function invariants associated to complex simple Lie group} \label{subsec:intro_L-func_inv_general_Lie}

% --------------------------------------------------------------------------

Next, we consider the $ G $ WRT invariants $ Z_k^G(M) $ for general complex simple Lie group $ G $.
In this case, Park~\cite{Park} introduced the extension of GPPV invariants for negative definite plumbed manifolds $ M $.
Chung~\cite{Chung:SU(N)} also introduced the extension of GPPV invariants for Seifert homology spheres $ M $ by a different expression.
Murakami--Terashima~\cite{M-Terashima} extended Chung's construction for general complex simple Lie groups $ G $ and Seifert homology spheres $ M $ as partial theta functions for a degenerate quadratic form.
Murakami--Terashima's $ q $-series are expected to coincide with Park's one.
Murakami--Terashima~\cite{M-Terashima} proved that its radial limits converge and can be written as finite sums which coincide with WRT invariants when $ \frakg $ is simply-laced.
We will give precise statements later in \cref{subsec:general_Lie}.

In this case, we introduce the $ L $-function invariant as
\[
L_b^{\frakg} (s; M, k)
\coloneqq
\frac{1}{\Gamma(s)} 
\int_{0}^{\infty} \Phi \left( \zeta_k e^{-t}\right) t^{s-1} dt,
\]
where $ \Phi(q) $ is Murakami--Terashima's $ q $-series.
We will express $ \Phi(q) $ and $ L_b^{\frakg} (s; M, k) $ explicitly later in \cref{subsec:general_Lie,subsec:proof_WRT_general_Lie} respectively.
We obtain the following statement as in \cref{thm:L-func_WRT}.

\begin{thm} \label{thm:L-func_WRT_general_Lie}
	For a general complex simple Lie algebra $ \frakg $ and a Seifert homology sphere $ M $, the $ L $-function $ L_b^{\frakg} (s; M, k) $ extends holomorphically on $ \bbC $.
	Moreover, if $ \frakg $ is simply-laced, then we have 
	\begin{align}
		Z_{k}^{\frakg}(M) 
		=
		&\frac{ (-1)^{\abs{\Delta_+}} \sqrt{[X \colon Y]} \zeta_8^{\dim_{\bbC} \frakg}}{\abs{W}}
		\left( \sqrt{k} \zeta_8^{-1} \right)^{\dim_{\bbC} \frakh}
		\lim_{q \to \zeta_k} L_b^{\frakg} (0; M, k).
	\end{align}
	Here, we use notations which we prepare in \cref{subsec:general_Lie}.
\end{thm}

% --------------------------------------------------------------------------

\subsection{Linear relations between special values at negative integers of some $ L $-functions} \label{subsec:intro_Bernoulli}

% --------------------------------------------------------------------------

We can prove \cref{thm:L-func_WRT,thm:L-func_WRT_general_Lie} by combining these two techniques to compute asymptotic expansions: the Euler--Maclaurin summation formula and Mellin transform.
We can apply this method for a general setting independent of the WRT invariants.

Let $ N $ be a positive integer, $ \lambda \in \Z^N $ and $ F \colon \Z_{\ge 0}^N + \lambda \to \bbC $ be an element of 
\[ 	
\delta(y \in \Z_{\ge 0}^N + \lambda) \cdot \bbC[ y_i, \delta(y_i \in a + k \Z) \mid 1 \le i \le N, a, k \in \Z ],	
\]
where $ \delta $ be the Kronecker delta function and we stipulate that $ a + k\Z = \{ a \} $ when $ k=0 $.	
Let $ \alpha \in \R_{> 0}^N $ be a vector,
$ w \colon \R^N \to \R_{\ge 0} $ be a $ C^\infty $ function such that $ w(n + \alpha) \neq 0 $ for $ n \in \Z^N $ 
and $ e^{-w(x)} $ is of rapid decay as $ x_i \to \infty $ for each $ 1 \le i \le N $, that is, $ x_i^m \frac{\partial^n f}{\partial x_i^n} (x) $ is bounded as $ x_i \to \infty $ for any $ m, n \in \Z_{\ge 0} $.
We define two $ L $-functions as
\begin{align}
	L(s_1, \dots, s_N; F, \alpha)
	&\coloneqq 
	\sum_{l \in \lambda + \Z_{\ge 0}^N} \frac{F(l)}{(n_1 + \alpha_1)^{s_1} \cdots (n_N + \alpha_N)^{s_N}}, \\
	L(s; F, \alpha, w)
	&\coloneqq 
	\sum_{l \in \Z_{\ge 0}^N} \frac{F(l)}{w(n + \alpha)^{s}}
\end{align}
for complex variables $ s, s_1, \dots, s_N $ with sufficiently large real parts.
Here, we remark that Hurwitz zeta functions and Barnes zeta functions are special cases of $ L(s_1, \dots, s_N; F, \alpha) $ and Epstein zeta functions can be written as finite sums of special cases of $ L(s; F, \alpha, w) $.

Our results are the following statement.

\begin{thm} \label{thm:relation_Bernoulli}
	The $ L $-functions $ L_{C, P, \alpha}(s_1, \dots, s_N) $ and $ L(s; C, P, \alpha, w) $ extend meromorphic functions in each variable on $ \bbC $ and  for each $ M \in \Z $ it holds that
	\begin{align}
		\Res_{s=-M} \Gamma(s) L(s; F, \alpha, w)
		=
		\sum_{m \in \Z^N, \, m_1 + \cdots + m_N = M}
		&(-1)^{m_1 + \cdots + m_N} 
		\frac{\partial^M e^{-w(x_1, \dots, x_N)}}{\partial x_1^{m_1} \cdots \partial x_N^{m_N}} (0, \dots, 0)
		\\
		&\cdot \Res_{s_1 = -m_1} \cdots \Res_{s_N = -m_N} 
		\Gamma(s_1) \cdots \Gamma(s_N)
		L(s_1, \dots, s_N; F, \alpha).
	\end{align}
	Here, for a $ C^\infty $ function $ g \colon \R \to \bbC $ of rapid decay as $ x \to \infty $, define
	\[
	%g^{(-1)}(x) = 
	\frac{d^{-1}g}{d x^{-1}}(x) \coloneqq -\int_x^\infty g(x') dx'.
	\]
\end{thm}

We can express values of these $ L $-functions at negative integers as finite sums of Bernoulli polynomials as described in \cref{rem:generalized_Bernoulli} later.

% --------------------------------------------------------------------------

\subsection{Organization of this paper} \label{subsec:Organization}

% --------------------------------------------------------------------------

In \cref{sec:asymp}, we recall two techniques to compute asymptotic expansions: the Euler--Maclaurin summation formula and Mellin transform.
In \cref{sec:relation_Bernoulli}, we prove \cref{thm:relation_Bernoulli}.
In \cref{sec:homological_blocks}, we summarize the GPPV invariants and Murakami--Terashima's $ q $-series.
In \cref{sec:proof}, we prove \cref{thm:L-func_WRT,thm:L-func_WRT_general_Lie}.

% --------------------------------------------------------------------------

\section*{Acknowledgement} \label{sec:acknowledgement}

% --------------------------------------------------------------------------

The author would like to show great appreciation to Masanobu Kaneko and Takuya Yamauchi for giving so much advice. 
The author would like to thank Yuji Terashima, Kazuhiro Hikami, Toshiki Matsusaka and William Misteg\aa{}rd for giving many comments. 
The author is supported by JSPS KAKENHI Grant Number JP23KJ1675.

% --------------------------------------------------------------------------

\section{Asymptotic formulas} \label{sec:asymp}

% --------------------------------------------------------------------------

In this section, we recall asymptotic formulas.

% --------------------------------------------------------------------------

\subsection{An asymptotic formula derived from Euler--Maclaurin summation formula} \label{subsec:asymp_EM}

% --------------------------------------------------------------------------

To begin with, we prepare the notation for asymptotic expansion by Poincar\'{e}. 

\begin{dfn}[Poincar\'{e}]
	Let $ L $ be a positive number, $ f \colon \R_{>0} \to \bbC $ be a map, $ t $ be a variable of $ \R_{>0} $, and $ (a_n)_{n \ge -L} $ be a family of complex numbers.
	Then, we write
	\[
	f(t) \sim \sum_{n \ge -L} a_n t^{n} \text{ as } t \to +0
	\]
	if for any positive number $ M $ there exist positive numbers $ K_M $ and $ \varepsilon $ such that
	\[
	\abs{ f(t) - \sum_{-L \le n \le M} a_n t^{n} }
	\le K_M \abs{ t^{M+1}}
	\]
	for any $ 0 < t < \varepsilon $.
	In this case, we call the infinite series $ \sum_{n \ge -L} a_n t^{n} $ as the \textbf{asymptotic expansion} of $ f(t) $ as $ t \to +0 $. 
\end{dfn}

Here, we remark that asymptotic expansions are typically divergent series.

We also need the following terminologies.

\begin{dfn}
	A $ C^\infty $ function $ f \colon \R \to \bbC $ is called of \textbf{rapid decay} as $ x \to \infty $ if $ x^m f^{(n)} (x) $ is bounded as $ x \to \infty $ for any $ m, n \in \Z_{\ge 0} $.
\end{dfn}

For functions $ g \colon \R \to \bbC $ and $ f \colon \R^N \to \bbC $, we define
\[
g^{(-1)}(x) = \frac{d^{-1}}{d x^{-1}} g(x) \coloneqq -\int_x^\infty g(x') dx', \quad
f^{(n)}(x) \coloneqq \frac{\partial^{n_1 + \cdots + n_N} f}{\partial x_1^{n_1} \cdots \partial x_r^{n_N}} (x).
\]

\begin{dfn} \label{dfn:Hadamard} 
	Let $ N $ be a positive integer. For a formal Laurent series	
	\[	
	\varphi(t_1, \dots, t_N) = \sum_{m \in \Z^N} B_m t_1^{m_1} \dots t_N^{m_N} \in \bbC((t_1, \dots, t_N))	
	\] 
	and a $ C^\infty $ function $ f \colon \R^{N} \to \bbC $ such that $ f^{(m)} (0) $ converges for any $ m \in \Z^N $, define their Hadamard product $ \varphi \odot f (t) \in \bbC((t)) $ as 
	\[	
	\varphi \odot f (t)	\coloneqq	\sum_{m \in \Z^N} B_m f^{(m)}(0) t^{m_1 + s + m_N}.	
	\]
\end{dfn}

Our key asymptotic formula is the following statement.

\begin{prop}[{\cite[Proposition 3.6]{M:GPPV}}] \label{prop:asymp_lim} 
	Let $ N $ be a positive integer, $ \lambda \in \Z^N $ and $ F \colon \Z_{\ge 0}^N + \lambda \to \bbC $ be an element of 
	\[ 	
	\delta(y \in \Z_{\ge 0}^N + \lambda) \cdot \bbC[ y_i, \delta(y_i \in a + k \Z) \mid 1 \le i \le N, a, k \in \Z ],	
	\]
	where $ \delta $ be the Kronecker delta function and we stipulate that $ a + k\Z = \{ a \} $ when $ k=0 $.
	Let 
	\[	
	\varphi_{F, u} (t_1, \dots, t_N)	\coloneqq	\sum_{l \in \Z_{\ge 0}^N + \lambda} F(l) 	e^{t_1 (l_1 + u_1) + s + t_N (l_N + u_N)}	\in \bbC ((t_1, \dots, t_N))[u_1, \dots, u_N].
	\]	
	Then, for any $ \alpha \in \R^N $ and any $ C^\infty $ function $ f \colon \R^{N} \to \bbC $ of rapid decay as $ x_1, \dots, x_N \to \infty $,	the following asymptotic formula holds as $t$ tends to $+0$	
	\[	
	\sum_{l \in \Z_{\ge 0}^N + \lambda} F(l) f(t(l+\alpha))	
	\sim \varphi_{F, \alpha} \odot f (t).
	\]
\end{prop}

\cref{prop:asymp_lim} generalizes \cite[Equation (44)]{Zagier:asymptotic}, \cite[Equation (2.8)]{BKM} and \cite[Lemma 2.2]{BMM:high_depth}.\cref{prop:asymp_lim} was proved by using the Euler--Maclaurin summation formula.
Moreover, as shown in the following, \cref{prop:asymp_lim} is a generalization of the following well-known two summation formulas.

\begin{rem}
	\begin{enumerate}
		\item In the case when $ N=1 $ and $ F(y) = \delta(y \in \Z_{\ge 0}) $, since
		\[
		\varphi_{F, \alpha} (t)
		=
		\sum_{l=0}^\infty e^{t(l+\alpha)}
		= 
		\frac{e^{\alpha t}}{1 - e^t}
		=
		\sum_{m=-1}^\infty \frac{B_{m+1} (\alpha)}{(m+1)!} t^m,
		\]
		where $ B_m(\alpha) $ is the $ m $-th Bernoulli polynomial,
		\cref{prop:asymp_lim} claims
		\begin{equation} \label{eq:EM_1variable}
			\sum_{n = 0}^{\infty} f(t(n+\alpha))
			\sim \sum_{m = -1}^{\infty} \frac{B_{m+1}(\alpha)}{(m+1)!}
			f^{(m)}(0) t^{m},
		\end{equation}
		which follows by the Euler--Maclaurin summation formula.
		
		\item In the case when $ N=1 $ and $ F(y) = (-1)^n \delta(y \in \Z_{\ge 0}) $, since
		\[
		\varphi_{F, \alpha} (t)
		=
		\sum_{l=0}^\infty (-1)^n e^{t(l+\alpha)}
		= 
		\frac{e^{\alpha t}}{1 + e^t}
		=
		\frac{1}{2} \sum_{m=0}^\infty \frac{E_{m} (\alpha)}{m!} t^m,
		\]
		where $ E_m(\alpha) $ is the $ m $-th Euler polynomial,
		\cref{prop:asymp_lim} claims
		\[
		\sum_{n = 0}^{\infty} (-1)^n f(t(n+\alpha))
		\sim 
		\frac{1}{2} \sum_{m = 0}^{\infty} \frac{E_{m}(\alpha)}{m!}
		f^{(m)}(0) t^{m},
		\]
		which follows by the Euler--Boole summation formula.
	\end{enumerate}
\end{rem}

\begin{rem} \label{rem:generalized_Bernoulli}
	We can regard the Laurent coefficients of $ \varphi_{F, \alpha} (t) $ in \cref{prop:asymp_lim} at $ t=0 $ as the finite sum of Bernoulli polynomials in $ \alpha_1, \dots, \alpha_N $.
	Thus, we can consider them as a generalization of Bernoulli polynomials.
\end{rem}

% --------------------------------------------------------------------------

\subsection{An asymptotic formula derived from Mellin transform} \label{subsec:asymp_Mellin}

% --------------------------------------------------------------------------

The following proposition is proved by using the Mellin transform.

\begin{prop}[{Lawrence--Zagier~\cite[pp.~98, Proposition]{Lawrence-Zagier}, Zagier~\cite[Section 1]{Zagier:asymptotic}}] \label{prop:L-func_asymp}
	Let $ \Lambda $ be a set and
	$ (a_\lambda)_{ \lambda \in \Lambda } $ and $ (b_\lambda)_{ \lambda \in \Lambda } $ be sequences of complex numbers and positive real numbers respectively.
	Assume that the series
	\[
	\varphi(t) \coloneqq \sum_{\lambda \in \Lambda} a_\lambda e^{-b_\lambda t}
	\]
	converges uniformly on any compact subsets of $ \R_{> 0} $ and admits asymptotic expansions
	\[
	\varphi(t) \sim \sum_{m = -M}^{\infty} c_m t^m
	\]
	as $ t \to +0 $ for some  $ M \in \Z $ and it holds that $ \varphi(t) = O(e^{-ct}) $ as $ t \to \infty $ for some $ c > 0 $.
	Then, the $ L $-function
	\[
	L(s; \varphi) \coloneqq \sum_{\lambda \in \Lambda} \frac{a_\lambda}{{b_\lambda}^s}
	\]
	extends to the meromorphic function on $ \bbC $ and its possible poles lies on $ \Z_{\le M} $.
	Moreover, for any $ n \in \Z $, it holds that
	\begin{align}
		c_n = \Res_{s=-n} \Gamma(s) L(s; \varphi)
		=
		\begin{dcases}
			\frac{(-1)^n}{n!} L(-n; \varphi), & n \ge 0, \\
			(-n-1)! \Res_{s=-n} L(s; \varphi), & n < 0.
		\end{dcases}
	\end{align}
\end{prop}

% --------------------------------------------------------------------------

\section{Relations between generalized Bernoulli polynomials} \label{sec:relation_Bernoulli}

% --------------------------------------------------------------------------

% --------------------------------------------------------------------------

\subsection{Proof of \cref{thm:relation_Bernoulli}} \label{subsec:proof_relation_Bernoulli}

% --------------------------------------------------------------------------

By comparing with \cref{prop:asymp_lim,prop:L-func_asymp}, we can prove \cref{thm:relation_Bernoulli}.

\begin{proof}[Proof of \cref{thm:relation_Bernoulli}]
	We use notations in \cref{thm:relation_Bernoulli}.
	Let $ f(x_1, \dots, x_N) \coloneqq e^{-w(x_1, \dots, x_N)} $ and
	\[
	\sum_{m \in \Z^{N}} B_m t_1^{m_1} \cdots t_{N}^{m_{N}}
	\coloneqq
	\varphi_{F, \alpha} (t_1, \dots, t_N).
	\]
	By comparing with \cref{prop:asymp_lim,prop:L-func_asymp}, for any $ M \in \Z $ we have
	\[
	\Res_{s=-M} \Gamma(s) L(s; F, \alpha, w)
	=
	\sum_{m \in \Z^N, \, m_1 + \cdots + m_N = M}
	(-1)^{m_1 + \cdots + m_N} 
	\frac{\partial^M e^{-w(x_1, \dots, x_N)}}{\partial x_1^{m_1} \cdots \partial x_N^{m_N}} (0, \dots, 0)
	B_m.
	\]
	By using \cref{prop:asymp_lim,prop:L-func_asymp} inductively on $ N $, we have
	\[
	B_m = 
	\Res_{s_1 = -m_1} \cdots \Res_{s_N = -m_N} 
	\Gamma(s_1) \cdots \Gamma(s_N)
	L(s_1, \dots, s_N; F, \alpha).
	\]
	Thus, we obtain the claim.
\end{proof}

% --------------------------------------------------------------------------

\subsection{Examples} \label{subsec:ex_relation_Bernoulli}

% --------------------------------------------------------------------------

\begin{ex}[Partial Epstein zeta functions]
	Let $ Q \colon \R^N \to \R_{\ge 0} $ be a positive definite quadratic form and $ \alpha \in \R^N_{>0} $.
	Define a partial Epstein zeta function for $ Q $ as 
	\[
	L(s; \alpha, Q)
	\coloneqq
	\sum_{n \in \Z_{\ge 0}^N} \frac{1}{Q(n + \alpha)^{s}}.
	\]
	Then, for any $ M \in \Z $ we have
	\begin{align}
		&\Res_{s=-M} \Gamma(s) L(s; \alpha, Q)
		\\
		= \,
		&\sum_{m \in \Z_{\ge -1}^N, \, m_1 + \cdots + m_N = M}
		(-1)^{m_1 + \cdots + m_N} 
		\frac{\partial^m e^{-Q(x_1, \dots, x_N)}}{\partial x_1^{m_1} \cdots \partial x_N^{m_N}} (0, \dots, 0)
		\frac{B_{m_1+1}(\alpha_1)}{(m_1 + 1)!} \cdots \frac{B_{m_N+1}(\alpha_N)}{(m_N + 1)!}.
	\end{align}
\end{ex}

\begin{ex}[Partial zeta functions attached to binary cubic forms]
	Let $ a, b, c, d \in \R_{>0} $ and $ \alpha \in \R^N_{>0} $.
	Define a binary cubic form $ w \colon \R^2 \to \R $ as $ w(x, y) \coloneqq ax^3 + bx^2 y + cxy^2 + dy^3 $ and its partial zeta function as 
	\[
	L(s; \alpha, w)
	\coloneqq
	\sum_{n \in \Z_{\ge 0}^2} \frac{1}{w(n + \alpha)^{s}}.
	\]
	Then, for any $ M \in \Z $ we have
	\begin{align}
		&\Res_{s=-M} \Gamma(s) L(s; \alpha, w)
		\\
		= \,
		&\sum_{m \in \Z_{\ge -1}^2, \, m_1 + m_2 = M}
		(-1)^{m_1 + m_2} 
		\frac{\partial^m e^{-w(x_1, x_2)}}{\partial x_1^{m_1} \partial x_2^{m_2}} (0, 0)
		\frac{B_{m_1+1}(\alpha_1)}{(m_1 + 1)!} \frac{B_{m_2+1}(\alpha_2)}{(m_2 + 1)!}.
	\end{align}
\end{ex}

% --------------------------------------------------------------------------

\section{Homological blocks} \label{sec:homological_blocks}

% --------------------------------------------------------------------------

In this section, we recall two cases when radial limit conjecture is solved.

% --------------------------------------------------------------------------

\subsection{Gukov--Pei--Putrov--Vafa invariants} \label{subsec:GPPV_inv}

% --------------------------------------------------------------------------

Let $ \Gamma = (\Vert, E, (b_I)_{I \in \Vert}) $ be a plumbing graph, that is, a finite tree with the vertex set $ \Vert $, the edge set $ E $, and integral weights $ b_I \in \Z $ for each vertex $ I \in \Vert $.
Then, we obtain a framed link $ \calL(\Gamma) $ whose all components are trivial knots corresponding to each vertex $ I $ with framing $ b_I $ and two components corresponding to $ I $ and $ I' $ chained together if and only if there exists an edge between $ I $ and $ I' $.
We show an example in \cref{fig:surgery_diagram}.

\begin{figure}[tb]
	\begin{minipage}[b]{0.45\linewidth}
		\centering
		\begin{tikzpicture}
			\node[shape=circle,fill=black, scale = 0.4] (1) at (0,0) { };
			\node[shape=circle,fill=black, scale = 0.4] (2) at (1.5,0) { };
			\node[shape=circle,fill=black, scale = 0.4] (3) at (-1,-1) { };
			\node[shape=circle,fill=black, scale = 0.4] (4) at (-1,1) { };
			\node[shape=circle,fill=black, scale = 0.4] (5) at (2.5,1) { };
			\node[shape=circle,fill=black, scale = 0.4] (6) at (2.5,-1) { };
			
			\node[draw=none] (B1) at (0,0.4) {$ b_1 $};
			\node[draw=none] (B2) at (1.5, 0.4) {$ b_2 $};
			\node[draw=none] (B3) at (-0.6,1) {$ b_3 $};
			\node[draw=none] (B4) at (-0.6,-1) {$ b_4 $};
			\node[draw=none] (B5) at (2.1,1) {$ b_5 $};		
			\node[draw=none] (B6) at (2.1,-1) {$ b_6 $};	
			
			\path [-](1) edge node[left] {} (2);
			\path [-](1) edge node[left] {} (3);
			\path [-](1) edge node[left] {} (4);
			\path [-](2) edge node[left] {} (5);
			\path [-](2) edge node[left] {} (6);
		\end{tikzpicture}
	\end{minipage}
	\begin{minipage}[b]{0.45\linewidth}
		\centering
		\begin{tikzpicture}[scale=0.5]
			\begin{knot}[
				clip width=5, 
				flip crossing=2, 
				flip crossing=4, 
				flip crossing=6, 
				flip crossing=7, 
				flip crossing=9
				]
				\strand[thick] (0, 0) circle [x radius=3cm, y radius=1.5cm];
				\strand[thick] (0, 0) +(4cm, 0pt) circle [x radius=3cm, y radius=1.5cm];
				\strand[thick] (0, 0) +(-2.5cm, 2cm) circle [radius=1.5cm];
				\strand[thick] (0, 0) +(-2.5cm, -2cm) circle [radius=1.5cm];
				\strand[thick] (0, 0) +(6.5cm, 2cm) circle [radius=1.5cm];
				\strand[thick] (0, 0) +(6.5cm, -2cm) circle [radius=1.5cm];
				
				\node (1) at (0, 2) {$b_{1}$};
				\node (2) at (4, 2) {$b_{2}$};
				\node (3) at (-2.5, 4) {$b_{3}$};
				\node (4) at (-2.5, -4) {$b_{4}$};
				\node (5) at (6.5, 4) {$b_{5}$};
				\node (6) at (6.5, -4) {$b_{6}$};
			\end{knot}
		\end{tikzpicture}
	\end{minipage}
	\caption{An example of plumbing graph $\Gamma$ and its corresponded framed link $ \calL(\Gamma) $.} 
	\label{fig:surgery_diagram}
\end{figure}
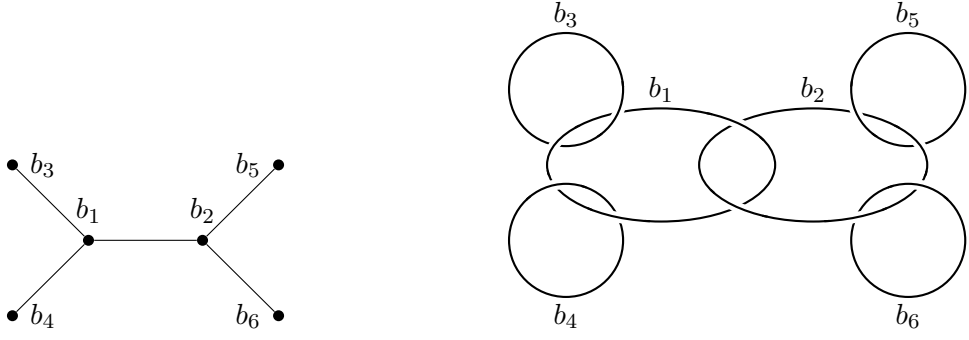

Let $ M $ be the 3-manifold obtained from $ S^3 $ through the surgery along the diagram $ \calL(\Gamma) $.
Such a manifold is called a plumbed manifold.
Following \cite{GPPV}, we assume that the linking matrix $ B $ of $ \Gamma $ is negative definite.
In this case, $ M $ is called a negative definite plumbed manifold.

Under the above setting, Gukov--Pei--Putrov--Vafa invariants are defined as follows.

\begin{dfn}[{\cite[Equation (A.29)]{GPPV}}]
	For a negative definite plumbed manifold $ M $, its Gukov--Pei--Putrov--Vafa invariant is defined as
	\[
	\widehat{Z}_{b} (q; M) 
	\coloneqq
	q^{-(3\abs{\Vert} + \tr B)/4} \,
	\mathrm{v.p.} \int_{\abs{x_{I}}=1, I \in \Vert} \Theta_{-B, b}(q; x)
	F\left( (x_I)_{I \in \Vert} \right) \prod_{I \in \Vert} \frac{dx_I}{2\pi\iu x_I}
	\]
	for $ b \in \Spin^c(M) \cong ( \delta + 2\Z^\Vert) /2B(\Z^\Vert) $, where $ \delta \coloneqq (\deg I)_{I \in \Vert} $,
	\[
	\Theta_{-B, b} (q; x) \coloneqq 
	\sum_{\ell \in b + 2W(\Z^\Vert)} q^{-\transpose{\ell} B \ell/4} \prod_{v \in V} x_I^{\ell_I}
	\]
	is the theta function defined for complex numbers $ q $ with $ \abs{q} < 1 $ and $ x  = (x_I)_{I \in \Vert} $ and
	\[
	F(x) = F\left( (x_I)_{I \in \Vert} \right) \coloneqq \prod_{I \in \Vert} F_I(x_I), \quad
	F_I(x_I) \coloneqq \left( x_I - x_I^{-1} \right)^{2 - \deg I}.
	\]
\end{dfn}

We express the GPPV invariants as false theta functions in the following.
For $ \ell = (\ell_I)_{I \in \Vert} \in \Z^\Vert $, define
\[
F_\ell \coloneqq \prod_{I \in \Vert} F_{I, \ell_I}, \quad
F_{I, \ell_I} \coloneqq \mathrm{v.p.} \int_{\abs{x_I} = 1} F_I(x_I) \frac{x_I^{\ell_I} dx_I}{2\pi\iu x_I},
\]
where
\[
\mathrm{v.p.} \int_{\abs{x} = 1} 
\coloneqq 
\frac{1}{2} \lim_{\veps \to +0} \left( \int_{\abs{x} = 1+\veps} + \int_{\abs{x} = 1-\veps} \right).
\]
By the calculation in \cite[p. 743]{Andersen-Mistegard}, we have 
\begin{equation} \label{eq:F_I_coeff_expression}
	F_{I, \ell_I} =
	\begin{dcases}
		-\ell_I, & \text{ if } \deg I = 1, \, \ell_I \in \{ \pm 1 \}, \\
		1, & \text{ if } \deg I = 2, \, \ell_I =0, \\
		\frac{1}{2} \binom{m + \deg I - 3}{\deg I - 3}, & \text{ if } \deg I \ge 3, \, \ell_I = \deg I -2 + 2m \text{ for some } m \in \Z_{\ge 0}, \\
		\frac{(-1)^{\deg I}}{2} \binom{m + \deg I - 3}{\deg I - 3}, & \text{ if } \deg I \ge 3, \, -\ell_I = \deg I -2 + 2m \text{ for some } m \in \Z_{\ge 0}, \\
		0, & \text{ otherwise}.
	\end{dcases}
\end{equation}

\begin{rem} \label{rem:F_l_symmetry}
	For any $ I \in \Vert $, we have $ F_{I, \ell_I} = (-1)^{\deg I} F_{I, \ell_I} $.
	Thus, we have 
	\begin{align}
		F_{-\ell}
		=
		(-1)^{\sum_{I \in \Vert} \deg I} F_\ell
		=
		F_\ell
	\end{align}
	by the handshaking lemma.
\end{rem}

\begin{rem} \label{rem:F_l_expansion}
	For each $ I \in \Vert $, we have
	\[
	\sum_{\ell_I \in \deg I + 2 \Z_{\ge -1}} F_{I, \ell_I} x_I^{\ell_I}
	=
	(-1)^{\deg I} F_I(x_I) \cdot
	\begin{dcases}
		1, \text{ if } \deg I \le 2, \\
		\frac{1}{2}, \text{ if } \deg I \ge 3
	\end{dcases}
	\]
	by \cref{eq:F_I_coeff_expression} and the binomial theorem.
	Thus, by the shakehand lemma, we have
	\[
	F\left( (x_I)_{I \in \Vert} \right)
	=
	2^{\abs{\Vert_{\ge 3}}}
	\sum_{\ell_I \in (\deg I)_{I \in \Vert} + 2 \Z_{\ge -1}^{\Vert}} F_{\ell} \prod_{I \in \Vert} x_I^{\ell_I},
	\]
	where $ \Vert_{\ge 3} \coloneqq \{ I \in \Vert \mid \deg I \ge 3 \} $.
\end{rem}

In the above notation, we can write
\begin{equation} \label{eq:GPPV_rep}
	\widehat{Z}_{b} (q; M) 
	\coloneqq
	q^{-(3\abs{\Vert} + \tr B)/4}
	\sum_{\ell \in b + 2B(\Z^\Vert)} F_\ell q^{-\transpose{\ell} B^{-1} \ell/4}.
\end{equation}

By considering GPPV invariants, Gukov--Pei--Putrov--Vafa~\cite{GPPV} refined radial limit conjecture (\cref{conj:radial_limit_su(2)}) for negative definite plumbed manifolds. 

\begin{conj}[{Radial limit conjecture, \cite[Equation (A.28)]{GPPV}}] \label{conj:GPPV}
	In the above situation, it holds that
	\[
	Z_k(M)
	=
	\sum_{a, b \in \Spin^c(M) / \{ \pm 1 \}} e^{2\pi\iu k \lk(a, a)} S_{a b}
	\lim_{q \to \zeta_k} \widehat{Z}_{b} (q; M),
	\]
	where 
	\[
	S_{a b} \coloneqq
	\frac{1}{2(\zeta_{2k} - \zeta_{2k}^{-1}) \sqrt{\abs{\det W}}}
	e^{-2\pi\iu \transpose{a} W^{-1} b}
	\]
\end{conj}

Here we remark that $ \Z^\Vert / B \Z^\Vert \cong H_1(M, \Z) $ and we can express the linking pairing $ \lk \colon H_1(M, \Z) \otimes H_1(M, \Z) \to \Q/\Z $ as
\[
\lk(a, b) = \transpose{a} B^{-1} b \bmod \Z.
\]

\cref{conj:GPPV} is proved by Murakami~\cite{M:GPPV}.

\begin{thm}[{Murakami~\cite{M:GPPV}}] \label{thm:GPPV_conj}
	$ \cref{conj:GPPV} $ is true.
\end{thm}

% --------------------------------------------------------------------------

\subsection{Homological blocks with simple Lie groups} \label{subsec:general_Lie}

% --------------------------------------------------------------------------

In this subsection, we present homological blocks with simple Lie groups $ G $ for Seifert fibered homology spheres $ M $, that is, a $ q $-series whose radial limits are $ G $ WRT invariants of $ M $.
As we mentioned in \cref{subsec:intro_L-func_inv}, such $ q $-series are constructed by Park~\cite{Park} for general $ G $ and negative definite plumbed manifolds $ M $, Chung~\cite{Chung:SU(N)} for $ G = \SU(N) $ and Seifert homology spheres $ M $ and Murakami--Terashima~\cite{M-Terashima} for general $ G $ and Seifert homology spheres $ M $.
Murakami--Terashima's homological block includes Chung's one and it is not checked that Park's homological block includes Murakami--Terashima's one, although it is expected.
Here, we recall the definition of homological blocks with simple Lie groups according to \cite{M-Terashima} since it is the only case when radial limit conjecture is solved.

To begin with, we prepare settings for Seifert fibered homology spheres.
Let $ n \ge 3 $ be an integer and $ (p_{1}, q_{1}), \dots, (p_{n}, q_{n}) $ be pairs of coprime integers such that $ 0 < q_1 < p_1, \dots, 0 < q_n < p_n $, $ p_1, \dots, p_n $ are pairwise coprime and
\[
P \sum_{1 \le i \le n} \frac{q_i}{p_i} = 1,
\]
where $ P \coloneqq p_1 \cdots p_n $.
We denote by $ M \coloneqq M(p_1/q_1, \dots, p_n/q_n) $ the Seifert fibered $ 3 $-manifold with $ n $-singular fibers and surgery integers $ p_1, \dots, p_n $ and $ q_1, \dots, q_n $.
This manifold is obtained by the surgery diagram in \cref{fig:surgery_diagram}.

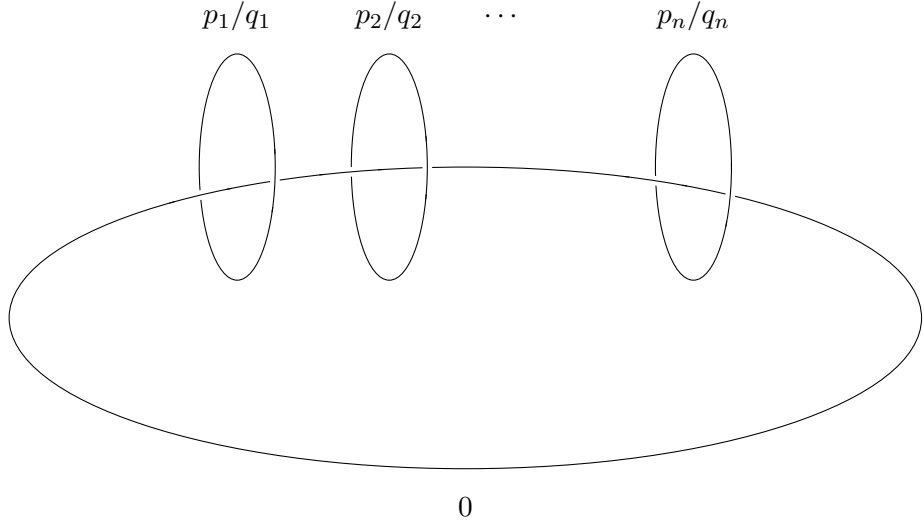
\begin{figure}[tb]
	\centering
	\begin{tikzpicture}
		\begin{knot}[
			%  draft mode=crossings,
			flip crossing/.list={2,4,6}
			]
			\strand (1,0) circle[x radius=6cm, y radius=2cm];
			\strand[] (-2,2) circle[x radius=0.5cm, y radius=1.5cm];
			\strand[] (0,2) circle[x radius=0.5cm, y radius=1.5cm];
			\strand[] (4,2) circle[x radius=0.5cm, y radius=1.5cm];
			
			\node (1) at (1, -2.5) {$ 0 $};
			\node (2) at (-2, 4) {$ p_1 / q_1 $};
			\node (3) at (0, 4) {$ p_2 / q_2 $};
			\node (4) at (1.5, 4) {$ \cdots $};
			\node (5) at (4, 4) {$ p_n / q_n $};
		\end{knot}
	\end{tikzpicture}
	\caption{The surgery diagram of $ M(p_1/q_1, \dots, p_n/q_n) $.} \label{fig:surgery_diagram}
\end{figure}

Let $ \frakg $ be a complex simple Lie algebra and $ \Delta_{+} $ be a set of positive roots.
Then, we define the homological block for $ M $ with $ \frakg $ as
\begin{equation} \label{eq:express_Phi}
	\Phi(q)
	\coloneqq
	q^{-(\dim_{\bbC} \frakg) \phi \abs{\rho}^{2}/2}
	\sum_{m \in \Z_{\ge m_0}^{\Delta_+}} q^{\abs{\sum_{\alpha \in \Delta_+} m_\alpha \alpha }^{2} / 8P}
	\prod_{\alpha \in \Delta_+} \chi_{m_\alpha},
\end{equation}
where 
\begin{itemize}
	\item denote
	\[
	\phi \coloneqq 3 - \frac{1}{P} + 12 \sum_{i=1}^n s(q_i, p_i),
	\]
	where 
	\[
	s(q_i, p_i) \coloneqq \frac{1}{4q_i} \sum_{j=1}^{q_i - 1}
	\cot \left( \frac{\pi j}{q_i} \right) \cot \left( \frac{\pi j p_i}{q_i} \right)
	\]
	are the Dedekind sums,
	\item denote $ \rho \coloneqq \sum_{\alpha \in \Delta_{+}} \alpha/2 \in \frakh_\R^* $, 
	\item define $ \abs{ \cdot } $ as follows:
	\begin{itemize}
		\item let $ \frakh $ be the Cartan subalgebra of $ \frakg $;
		\item let $ \frakh_\R^* $ be the real vector subspace of the dual space $ \frakh^* $ generated by the root system $ \Delta $;
		\item let $ \sprod{\cdot, \cdot} $ be the standard inner product on $ \frakh_\R^* $ normalized such that the longest roots have length $ \sqrt{2} $;
		\item for $ x \in \frakh_\R^* $, denote $ \abs{x} \coloneqq \sqrt{\sprod{x, x}} $;
	\end{itemize}	
	\item denote
	\[
	\sum_{m = m_0}^\infty \chi_m q^{m/2P}
	\coloneqq
	G_{p_1, \dots, p_n} (q)
	\coloneqq
	(q^{1/2} - q^{-1/2})^{2-n} \prod_{1 \le i \le n} (q_\alpha^{1/2p_i} - q_\alpha^{-1/2p_i}).
	\]
\end{itemize}

Here, we remark that we can express $m_0$ and $\chi_m$ as
\begin{align}
	m_0 &= P \left( - \frac{1}{p_1} - \cdots - \frac{1}{p_n} + n - 2 \right), \\
	\chi_m &=
	\begin{dcases}
		(-1)^n \veps_1 \cdots \veps_n \binom{m'+n-3}{n-3},
		& 
		\begin{gathered}
			\text{ if } \frac{m}{P} = \frac{\veps_1}{p_1} + \cdots + \frac{\veps_n}{p_n} + n - 2 + 2m'
			\text{ for some } \\
			\veps_1, \dots, \veps_n \in \{ \pm 1 \} \text{ and } m' \in \Z_{\ge 0},
		\end{gathered}
		\\
		0 & \text{ otherwise}.
	\end{dcases}
	\label{eq:express_chi}
\end{align}
by \cite[Remark 2.1]{M-Terashima} and the homological block $\Phi(q)$ is independent of a choice of a set of positive roots $ \Delta_{+} $ by \cite[Remark 2.2]{M-Terashima}.

Murakami--Terashima~\cite{M-Terashima} gave the following formula for radial limits of the homological block $\Phi(q)$, which solves radial limit conjecture for simply-laced $ \frakg $ and Seifert homology spheres $ M $.

\begin{thm} \label{thm:M-Terashima}
	Let $ k $ be a positive integer.
	\begin{enumerate}
		\item \textup{(\cite[Theorem 2.3]{M-Terashima})} It holds that
		\begin{align}
			\lim_{q \to \zeta_k} \Phi(q)
			= \,
			&\bm{e} \left( 
			-\frac{(\dim_{\bbC} \frakg) \phi \abs{\rho}^{2} }{2k}
			\right)
			\left( \frac{\zeta_8}{\sqrt{k}} \right)^{\dim_{\bbC} \frakh} 
			\frac{1}{\sqrt{[X : Y]}}
			\\
			&\sum_{\lambda \in X/k PY \smallsetminus \calM}
			\bm{e} \left( -\frac{1}{2Pk} \abs{\lambda}^{2} \right)
			\prod_{\alpha \in \Delta_+}
			G_{p_1, \dots, p_n} \left( \zeta_k^{\sprod{\lambda, \alpha}} \right),
		\end{align}
		where 
		\begin{itemize}
			\item denote $ \bm{e}(z) \coloneqq e^{2\pi\iu z} $ for a complex number $ z $,
			\item denote $ \zeta_l^{z} \coloneqq \bm{e}(z/l) $ for a positive integer $l$ and a complex number $z$, 
			%    \item 
			%    \[
			%	G_{p_1, \dots, p_n}(q)
			%	\coloneqq
			%	(q^{1/2} - q^{-1/2})^{2-n}
			%	\prod_{1 \le i \le n} (q^{1/2p_i} - q^{-1/2p_i}),
			%	\]
			\item $ X $ and $ Y $ are the weight and root lattices respectively, 
			\item $ \calM \coloneqq
			\{ \lambda \in X \mid 
			\sprod{\lambda, \alpha} \in k \Z \} $.
		\end{itemize}
		
		\item \textup{(\cite[Corollary 5.1]{M-Terashima})} Let $ Z_{k}^{\frakg}(M) \in \bbC $ be the WRT invariant of $ M $ with $ \frakg $.
		Then, when $\frakg$ is simply-laced, we have
		\begin{align}
			Z_{k}^{\frakg}(M) 
			=
			&\frac{ (-1)^{\abs{\Delta_+}} \sqrt{[X \colon Y]} \zeta_8^{\dim_{\bbC} \frakg}}{\abs{W}}
			\left( \sqrt{k} \zeta_8^{-1} \right)^{\dim_{\bbC} \frakh}
			\lim_{q \to \zeta_k} \Phi(q),
		\end{align}
		where $ W $ is the Weyl group.
	\end{enumerate}
\end{thm}

% --------------------------------------------------------------------------

\section{Proofs of main results} \label{sec:proof}

% --------------------------------------------------------------------------

% --------------------------------------------------------------------------

\subsection{Proof of \cref{thm:L-func_WRT,thm:L-func_CGP}} \label{subsec:proof_WRT}

% --------------------------------------------------------------------------

By \cref{eq:GPPV_rep}, we have
\begin{align}
	L_{b} (s; M, k)
	&\coloneqq
	\frac{1}{\Gamma(s)} 
	\int_{0}^{\infty} \widehat{Z}_b \left( \zeta_k e^{-t}; M \right) t^{s-1} dt \\
	&=
	4^s \sum_{\ell \in b + 2B(\Z^\Vert)} 
	\frac{F_\ell \zeta_{4k}^{-3 \abs{\Vert} - \tr B - \transpose{\ell} B^{-1} \ell}}{(-3 \abs{\Vert} - \tr B - \transpose{\ell} B^{-1} \ell)^s}.
\end{align}
By the representation of $ F_\ell $ in \cref{eq:F_I_coeff_expression}, $ L_{b} (s; M, k) $ is a finite linear combination of the $ L $-functions $ L(s; C, P, \alpha, w) $ in \cref{thm:relation_Bernoulli}.

\begin{proof}[Proof of \cref{thm:L-func_WRT,thm:L-func_CGP}]
	First, we prove \cref{thm:L-func_WRT}.
	By the above argument and \cref{prop:asymp_lim}, the function
	\[
	\varphi(t) \coloneqq 
	\widehat{Z}_{b} (\zeta_k e^{-t}; M)
	=
	\sum_{\ell \in b + 2B(\Z^\Vert)} F_\ell \zeta_{4k}^{-3 \abs{\Vert} - \tr B - \transpose{\ell} B^{-1} \ell}
	e^{(3 \abs{\Vert} + \tr B + \transpose{\ell} B^{-1} \ell) t/4}
	\]
	satisfies the assumption of \cref{prop:L-func_asymp}.
	Thus, $ L_{b} (s; M, k) $ extends meromorphically on $ \bbC $ by \cref{prop:L-func_asymp}.
	Moreover, \cref{thm:GPPV_conj} implies that asymptotic coefficients as $ t \to +0 $ of
	\[
	\sum_{a, b \in \Spin^c(M) / \{ \pm 1 \}} \bm{e} \left( k \lk(a, a) \right) S_{a b}
	\widehat{Z}_{b} (\zeta_k e^{-t}; M)	
	\]
	vanishes at negative powers and equals to $ (\zeta_{2k} - \zeta_{2k}^{-1}) \WRT_k(M) $ at the constant term.
	Thus, the $ L $-function
	\[
	\sum_{a, b \in \Spin^c(M) / \{ \pm 1 \}} \bm{e} \left( k \lk(a, a) \right) S_{a b} L_b (s; M, k)	
	\]
	is entire and its special value at $ s=0 $ is equal to $ (\zeta_{2k} - \zeta_{2k}^{-1}) \WRT_k(M) $ by \cref{prop:L-func_asymp}.
	
	\cref{thm:L-func_CGP} follows from the same argument.
\end{proof}

% --------------------------------------------------------------------------

\subsection{Proof of \cref{thm:functional_equation}} \label{subsec:proof_functional_equation}

% --------------------------------------------------------------------------

To begin with, we state the modular transform of GPPV invariants of Seifert homology spheres.

\begin{thm}[{Andersen--Misteg\aa{}rd~\cite[p.751]{Andersen-Mistegard}}] \label{thm:AM_CS}
	For a Seifert homology sphere $ M $, it holds that
	\[
	(-1)^n q^{\Delta} \widehat{Z}_0(M; q)
	=
	2 \sqrt{\frac{2P \iu}{\tau}} \sum_{\theta \in \CS_{\bbC}^* (M)} Z_\theta^* \left( M; -\frac{1}{\tau} \right)
	+ \frac{1}{\sqrt{\tau}} \calL^{\pi/2 - \varepsilon} \widehat{\varphi} (\tau),
	\]
	where 
	$ p_1, \dots, p_n $ are pairwise coprime positive integers such that $ M $ is the oriented Seifert fibered integral homology sphere $ \Sigma(p_1, \dots, p_n) $,
	$ q_1, \dots, q_n $ are positive integers such that $ \gcd(p_i, q_i) = 1 $ and $ 1 \le q_i \le p_i $ for each $ 1 \le i \le n $ and
	\[
	\sum_{1 \le i \le n} \frac{q_i}{p_i}
	\equiv -\frac{1}{p_1 \cdots p_n} \bmod \Z,
	\]
	\[
	\frac{p_i}{q_i}
	=
	k_{i, 1} - \cfrac{1}{k_{i, 2} - \cfrac{1}{\ddots - \cfrac{1}{k_{i, s_i}}}}
	\]
	is the continued fraction expansion with $ k_{i,j} \in \Z_{\ge 2} $,
	\begin{align}
		P &\coloneqq p_1 \cdots p_n,
		\\
		\Delta &\coloneqq
		\sum_{a}
		\\
		\CS_{\bbC}^* (M) 
		&\coloneqq
		\left\{ \text{non-zero Chern--Simons invariants of $ M $} \right\}
		\subset \Q/\Z,
		\\
		G(z) &\coloneqq
		\left( z^{P} - z^{-P} \right)^{2-n} \prod_{1 \le i \le n} \left( z^{P/p_i} - z^{-P/p_i} \right),
		\\
		Z_\theta^* (M; \tau)
		&\coloneqq
		\frac{1}{2}
		\sum_{\substack{
				m \in \Z, \\
				\theta \equiv -m^2/4P \bmod \Z
		}}
		\sgn(m) \Res_{x=m} q^{x^2/4P} G \left( e^{\pi\iu x/P} \right),
		\\
		\widehat{\varphi}(\xi)
		&\coloneqq
		\frac{1}{\sqrt{\pi \xi}} G \left( \exp \left( \sqrt{\frac{2\pi \iu \xi}{P}} \right) \right),
		\\
		\calL^{\pi/2 - \varepsilon} \widehat{\varphi} (\tau) 
		&\coloneqq
		\int_{0}^{e^{(\pi/2 - \varepsilon) \iu} \infty} e^{-\xi / \tau} \widehat{\varphi}(\xi) d \xi
	\end{align}
	and $ \varepsilon > 0 $ is a sufficiently small number.
\end{thm}

In \cite{Andersen-Mistegard}, \cref{thm:AM_CS} was used to give the asymptotic expansion of WRT invariants for Seifert homology spheres.

Here we remark that Andersen--Misteg\aa{}rd~\cite[Corollary 9]{Andersen-Mistegard} proved that
\begin{equation} \label{eq:CS_explicit}
	\CS_{\bbC}^* (M)
	= \left\{ - \frac{m^2}{4P} \bmod \Z \relmiddle| m \in \Z \text{ such that at most $ n-3 $ of $ 1 \le i \le n $ satisfies $ p_i \mid m $} \right\}.
\end{equation}

We describe the representation of the function $ Z_\theta^* \left( M; \tau \right) $.
Since $ G(z^{-1}) = G(z) $, for each $ m \in \Z $ we have
\[
\Res_{x=-m} q^{x^2/4P} G \left( e^{\pi\iu x/P} \right)
=
-\Res_{x=m} q^{x^2/4P} G \left( e^{\pi\iu x/P} \right).
\]
Thus, for each $ \theta \in \CS_{\bbC}^* (M) $, by taking $ 0 < l < P $ such that $ \theta \equiv -l^2/4P \bmod \Z $, we have
\begin{align}
	Z_\theta^* \left( M; \tau \right)
	&=
	\sum_{\substack{
			m \in \Z_{> 0}, \\
			\theta \equiv -m^2/4P \bmod \Z
	}}
	\Res_{x=m} q^{x^2/4P} G \left( e^{\pi\iu x/P} \right)
	\\
	&=
	\left(
	\sum_{m \in l + 2P\Z_{\ge 0}} + \sum_{m \in 2P - l + 2P\Z_{\ge 0}}
	\right)
	\Res_{x=m} q^{x^2/4P} G \left( e^{\pi\iu x/P} \right).
\end{align}
By the calculation in \cite[p.756]{Andersen-Mistegard}, there exists maps
\[
g_{0} (y), \dots, g_{n-2} (y) \in \bigoplus_{C \colon \Z/2P\Z \to \bbC} C(y) \cdot \bbC[y]
\]
such that for each $ \theta \in \CS_{\bbC}^* (M) $ it holds that
\begin{equation} \label{eq:Z*_rep}
	Z_\theta^* \left( M; \tau \right)
	=
	\sum_{d=0}^{n-2}
	\tau^d 
	\sum_{\substack{
			m \in \Z_{> 0}, \\
			\theta \equiv -m^2/4P \bmod \Z
	}}
	g_d (m) q^{m^2/4P}.
\end{equation}

\begin{proof}[Proof of \cref{thm:functional_equation}]
	%	We denote $ \bm{e}(z) \coloneqq e^{2\pi\iu z} $ for a complex number $ z $.
	By \cref{thm:AM_CS}, we have
	\begin{align}
		&\phant
		(-1)^n L_0 (s; M)
		\\
		&=
		\frac{\sqrt{-2\pi\iu}}{\Gamma(s)}
		\int_{0}^{\infty} e^{-\Delta t}
		\left(
		2 \sqrt{2P \iu} 
		\sum_{\theta \in \CS_{\bbC}^* (M)} Z_\theta^* \left( M; \frac{2\pi\iu}{t} \right)
		+ \calL^{\pi/2 - \varepsilon} \widehat{\varphi} \left( \frac{t\iu}{2\pi} \right)
		\right)
		t^{s-3/2} dt.
		%		\\
		%		&=
		%		\frac{\sqrt{-2\pi\iu}}{\Gamma(s)}
		%		\int_{0}^{\infty} e^{-\Delta t}
		%		\left(
		%			2 \sqrt{2P \iu} 
		%			\sum_{\theta \in \CS_{\bbC}^* (M)} 
		%			\sum_{\substack{
				%					m \in \Z, \\
				%					\theta \equiv -m^2/4P \bmod \Z
				%			}}
		%			\sgn(m) \Res_{x=m} \bm{e} \left( \frac{\pi\iu x^2}{2Pt} \right) G \left( \bm{e} \left( \frac{x}{2P} \right) \right)
		%		\right.
		%		\\
		%		&\phant
		%		\hphantom{
			%			\frac{\sqrt{-2\pi\iu}}{\Gamma(s)}
			%			\int_{0}^{\infty} e^{-\Delta t}
			%			\left( \right.
			%		}
		%		\left.
		%			\vphantom{\sum_{\substack{
					%						m \in \Z, \\
					%						\theta \equiv -m^2/4P \bmod \Z
					%			}}}
		%			+ \calL^{\pi/2 - \varepsilon} \widehat{\varphi} \left( \frac{t\iu}{2\pi} \right)
		%		\right)
		%		t^{s-3/2} dt
	\end{align}
	We have
	\begin{alignat}{2}
		&\phant
		\int_{0}^{\infty} e^{-\Delta t}
		\sum_{\theta \in \CS_{\bbC}^* (M)} Z_\theta^* \left( M; \frac{2\pi\iu}{t} \right)
		t^{s-3/2} dt
		& &
		\\
		&=
		\int_{0}^{\infty} e^{-\Delta t}
		\sum_{d=0}^{n-2}
		\left( \frac{t}{2\pi\iu} \right)^d 
		\sum_{m=1}^{\infty}
		g_d (m) e^{-(\pi m)^2 /Pt}
		t^{s-3/2} dt
		& &\quad \text{ by \cref{eq:Z*_rep}}
		\\
		&=
		\sum_{d=0}^{n-2}
		\left( 2\pi\iu \right)^{-d}
		\sum_{m=1}^{\infty}
		g_d (m) 
		\int_{0}^{\infty} e^{-\Delta t - (\pi m)^2 /Pt}
		t^{s+d-3/2} dt.
		& &
	\end{alignat}
	Since 
	\begin{equation} \label{eq:K-Bessel_integral_expression}
		K_s (z) = 
		\frac{1}{2} \left( \frac{2}{z} \right)^{s}
		\int_{0}^{\infty} e^{-t-z^2/4t} t^{s-1} dt,
		\quad
		\abs{z} < \frac{\pi}{4}
	\end{equation}
	by \cite[p.119, Equation (5.10.25)]{Lebedev}, this is equal to
	\[
	2\sum_{d=0}^{n-2}
	\left( 2\pi\iu \right)^{-d}
	\left( \frac{\pi}{\sqrt{\Delta P}} \right)^{s+d-1/2}
	\sum_{m=1}^{\infty}
	g_d (m) m^{s+d-1/2} K_{s+d-1/2} \left( 2\pi \sqrt{\frac{\Delta}{P}} m \right).
	\]
	This sum converges since $ K_{s+d-1/2} (x) = O(x^{-1/2} e^{-x}) $ as $ x \to \infty $ (\cite[Equation (5.11.9) in p.123]{Lebedev}).
	
	%	Here, for any $ m \in \Z $ with $ -m^2/4P \bmod \Z \notin \CS_{\bbC}^* (M) $, it holds that
	%	\[
	%	\Res_{x=m} q^{x^2/4P} G \left( \bm{e} \left( \frac{x}{2P} \right) \right)
	%	= 0
	%	\]
	%	since ...
	
	We also have
	\begin{align}
		&\phant
		\int_{0}^{\infty} e^{-\Delta t}
		\calL^{\pi/2 - \varepsilon} \widehat{\varphi} \left( \frac{t\iu}{2\pi} \right)
		t^{s-3/2} dt
		\\
		&=
		\int_{0}^{\infty} e^{-\Delta t}
		\left(
		\int_{0}^{e^{(\pi/2 - \varepsilon) \iu} \infty} e^{2\pi\iu\xi / t} 
		\frac{1}{\sqrt{\pi \xi}} G \left( \exp \left( \sqrt{\frac{2\pi \iu \xi}{P}} \right) \right)
		d \xi
		\right)
		t^{s-3/2} dt
		\\
		&=
		\int_{0}^{\infty} e^{-\Delta t}
		\left(
		\int_{0}^{e^{- \varepsilon \iu} \infty} e^{-x/t} 
		\sqrt{\frac{2}{x\iu}} G \left( \exp \left( \sqrt{-\frac{x}{P}} \right) \right)
		\frac{dx}{2\pi\iu}
		\right)
		t^{s-3/2} dt
		\\
		&=
		\frac{1}{\pi\sqrt{-2\iu}}
		\int_{0}^{e^{- \varepsilon \iu} \infty}
		\left(
		\int_{0}^{\infty} e^{-\Delta t - x/t}
		t^{s-3/2} dt
		\right)
		G \left( \exp \left( \sqrt{-\frac{x}{P}} \right) \right)
		\frac{dx}{\sqrt{x}}.
	\end{align}
	By \cref{eq:K-Bessel_integral_expression}, this is equal to
	\[
	\frac{\sqrt{2} \zeta_8^{5} \Delta^{-s/2 + 1/4}}{\pi}
	\int_{0}^{e^{- \varepsilon \iu} \infty}
	K_{s-1/2} (2\sqrt{\Delta x})
	G \left( \exp \left( \sqrt{-\frac{x}{P}} \right) \right)
	x^{s/2 - 3/4} dx.
	\]
	By substituting $ x $ with $ x^2 $, this is equal to
	\[
	\frac{2^{3/2} \zeta_8^{5} \Delta^{-s/2 + 1/4}}{\pi}
	\int_{0}^{e^{- \varepsilon \iu} \infty}
	K_{s-1/2} (2\sqrt{\Delta}x)
	G \left( \exp \left( -\frac{x\iu}{\sqrt{P}} \right) \right)
	x^{s - 1/2} dx.
	\]
	Finally, we prove this integral converges if $ \Re(s) > -1 $.
	Since $ K_\nu (x) = O(x^\nu + x^{-\nu}) $ and $ G(e^{x}) = O(x^2) $ as $ z \to 0 $, we have
	\[
	K_{s-1/2} (2\sqrt{\Delta}x)
	G \left( \exp \left( -\frac{x\iu}{\sqrt{P}} \right) \right)
	x^{s - 1/2}
	=
	\begin{dcases}
		O(x^2) & \text{ if } \Re(s) > \frac{1}{2}, \\
		O(x^{2s+1}) & \text{ if } \Re(s) < \frac{1}{2}.
	\end{dcases}
	\]
	Thus, the above integral converges if $ \Re(s) > -1 $.
\end{proof}

\begin{rem}
	In the right hand side in \cref{thm:functional_equation}, the integral
	\[
	\int_{0}^{e^{- \varepsilon \iu} \infty}
	K_{s-1/2} (2\sqrt{\Delta}x)
	G \left( \exp \left( -\frac{x\iu}{\sqrt{P}} \right) \right)
	x^{s - 1/2} dx
	\]
	cannot be expanded as an infinite sum
	\begin{equation} \label{eq:divergent_sum}
		\sum_{m=m_0}^{\infty} \chi_m I_m,
	\end{equation}
	here we define $ m_0 \in \Z $ and $ (\chi_m)_{m=m_0}^\infty $ as
	\[
	\sum_{m=m_0}^{\infty} \chi_m z^{m}
	\coloneqq
	G(z)
	=
	G(z^{-1})
	\]
	and
	\[
	I_m \coloneqq
	\int_{0}^{e^{- \varepsilon \iu} \infty}
	K_{s-1/2} (2\sqrt{\Delta}x)
	\exp \left( -m\frac{x\iu}{\sqrt{P}} \right)
	x^{s - 1/2} dx.
	\]
	In fact, the infinite sum \cref{eq:divergent_sum} diverges because
	\[
	\chi_m =
	\begin{dcases}
		\varepsilon_1 \cdots \varepsilon_n \binom{r+n-3}{r} & \text{ if } m= P \left( 2r+n-3 + \sum_{i=1}^n \varepsilon_i \frac{P}{p_i} \right)
		\text{ for some } \varepsilon_1, \dots, \varepsilon_n \in \{ \pm 1 \} \text{ and } r \in \Z_{\ge 0}, 
		\\
		0 & \text{ otherwise}
	\end{dcases}
	\]
	and $ I_m = O(m^{-1}) $ as $ m \to \infty $.
	Here, the last fact follows from $ K_\nu (x) = O(x^\nu + x^{-\nu}) $ as $ x \to 0 $ and generalized Watson's lemma (\cite[p.22]{Wong}), which states that
	\[
	\int_{0}^{e^{\theta \iu} \infty} f(x) e^{-zx} dx
	= O(z^{-\nu-1})
	\quad \text{ as } \abs{z} \to \infty
	\]
	for $ \theta \in \R $ and $ f \colon e^{\theta \iu} \R_{>0} \to \bbC $ with $ f(x) = O(x^\nu) $ as $ x \to 0 $.

\end{rem}

\subsection{Proof of \cref{thm:L-func_WRT_general_Lie}} \label{subsec:proof_WRT_general_Lie}

% --------------------------------------------------------------------------

By the expression of the homological block $ \Phi(q) $ in \cref{eq:express_Phi}, we have
\begin{align}
	L_b^{\frakg} (s; M, k)
	&\coloneqq
	\frac{1}{\Gamma(s)} 
	\int_{0}^{\infty} \Phi \left( \zeta_k e^{-t}\right) t^{s-1} dt \\
	&=
	\sum_{m \in \Z_{\ge m_0}^{\Delta_+}}
	\left( \frac{1}{8P} \abs{\sum_{\alpha \in \Delta_+} m_\alpha \alpha }^{2} - \frac{1}{2} (\dim_{\bbC} \frakg) \phi \abs{\rho}^{2} \right)^{-s}
	\prod_{\alpha \in \Delta_+} \chi_{m_\alpha}.
\end{align}
By the representation of $ \chi_{m} $ in \cref{eq:express_chi}, $ L_{b}^{\frakg} (s; M, k) $ is a finite linear combination of the $ L $-functions $ L(s; C, P, \alpha, w) $ in \cref{thm:relation_Bernoulli}.
Thus, we obtain a proof of \cref{subsec:proof_WRT} by the same argument in \cref{thm:L-func_WRT_general_Lie}.

% --------------------------------------------------------------------------
%		References
% --------------------------------------------------------------------------

\bibliographystyle{alpha}
\bibliography{quantum_inv,modular}

\newcommand{\etalchar}[1]{$^{#1}$}
\begin{thebibliography}{CGPM14}

\bibitem[AHJ{\etalchar{+}}17]{Andersen-Himpel-Jorgensen-Martens-McLellan}
J.~E. Andersen, B.~Himpel, S.~F. J\o{}rgensen, J.~Martens, and B.~McLellan.
\newblock The {W}itten-{R}eshetikhin-{T}uraev invariant for links in finite
  order mapping tori {I}.
\newblock {\em Adv. Math.}, 304:131--178, 2017.

\bibitem[AM22]{Andersen-Mistegard}
J.~E. Andersen and W.~Misteg{\aa}rd.
\newblock Resurgence analysis of quantum invariants of {S}eifert fibered
  homology spheres.
\newblock {\em Journal of the London Mathematical Society}, 105(2):709--764,
  2022.

\bibitem[BKM19]{BKM}
K.~Bringmann, J.~Kaszian, and A.~Milas.
\newblock Higher depth quantum modular forms, multiple {E}ichler integrals, and
  {$\mathfrak{sl}_3$} false theta functions.
\newblock {\em Res. Math. Sci.}, 6(2):Paper No. 20, 41, 2019.

\bibitem[BMM20]{BMM:high_depth}
K.~Bringmann, K.~Mahlburg, and A.~Milas.
\newblock Higher depth quantum modular forms and plumbed 3-manifolds.
\newblock {\em Lett. Math. Phys.}, 110(10):2675--2702, 2020.

\bibitem[CGP23]{Costantino-Gukov-Putrov}
F.~Costantino, S.~Gukov, and P.~Putrov.
\newblock Non-semisimple tqft's and bps q-series.
\newblock {\em SIGMA. Symmetry, Integrability and Geometry: Methods and
  Applications}, 19:010, 2023.

\bibitem[CGPM14]{CGP}
F.~Costantino, N.~Geer, and B.~Patureau-Mirand.
\newblock Quantum invariants of 3-manifolds via link surgery presentations and
  non-semi-simple categories.
\newblock {\em J. Topol.}, 7(4):1005--1053, 2014.

\bibitem[Chu22]{Chung:SU(N)}
H.-J. Chung.
\newblock {BPS} invariants for a knot in {S}eifert manifolds.
\newblock {\em Journal of High Energy Physics}, 2022(12):1--24, 2022.

\bibitem[FIMT21]{FIMT}
H.~Fuji, K.~Iwaki, H.~Murakami, and Y.~Terashima.
\newblock Witten–{R}eshetikhin–{T}uraev function for a knot in {S}eifert
  manifolds.
\newblock {\em Communications in Mathematical Physics}, 2021.

\bibitem[GM21]{Gukov-Manolescu}
S.~Gukov and C.~Manolescu.
\newblock A two-variable series for knot complements.
\newblock {\em Quantum Topology}, 12(1), 2021.
\newblock arXiv:1904.06057.

\bibitem[GPPV20]{GPPV}
S.~Gukov, D.~Pei, P.~Putrov, and C.~Vafa.
\newblock B{PS} spectra and 3-manifold invariants.
\newblock {\em J. Knot Theory Ramifications}, 29(2):2040003, 85, 2020.

\bibitem[GPV17]{GPV}
S.~Gukov, P.~Putrov, and C.~Vafa.
\newblock Fivebranes and $3$-manifold homology.
\newblock {\em J. High Energy Phys.}, 2017(071), 2017.

\bibitem[Hik05a]{Hikami:Bries}
K.~Hikami.
\newblock On the quantum invariant for the {B}rieskorn homology spheres.
\newblock {\em Internat. J. Math.}, 16(6):661--685, 2005.

\bibitem[Hik05b]{Hikami:Lattice}
K.~Hikami.
\newblock Quantum invariant, modular form, and lattice points.
\newblock {\em Int. Math. Res. Not.}, (3):121--154, 2005.

\bibitem[Hik06a]{Hikami:Seifert}
K.~Hikami.
\newblock On the quantum invariants for the spherical {S}eifert manifolds.
\newblock {\em Comm. Math. Phys.}, 268(2):285--319, 2006.

\bibitem[Hik06b]{Hikami:Lattice2}
K.~Hikami.
\newblock Quantum invariants, modular forms, and lattice points. {II}.
\newblock {\em J. Math. Phys.}, 47(10):102301, 32, 2006.

\bibitem[Hik11]{Hikami:radial_limit}
Decomposition of {W}itten-{R}eshetikhin-{T}uraev invariant: Author = {Hikami,
  Kazuhiro}, linking pairing and modular forms.
\newblock In {\em Chern-{S}imons gauge theory: 20 years after}, volume~50 of
  {\em AMS/IP Stud. Adv. Math.}, pages 131--151. Amer. Math. Soc., Providence,
  RI, 2011.

\bibitem[Leb72]{Lebedev}
N.~N. Lebedev.
\newblock {\em Special functions and their applications}.
\newblock Dover Publications, Inc., New York, 1972.
\newblock Revised edition, translated from the Russian and edited by Richard A.
  Silverman, Unabridged and corrected republication.

\bibitem[LZ99]{Lawrence-Zagier}
R.~Lawrence and D.~Zagier.
\newblock Modular forms and quantum invariants of {$3$}-manifolds.
\newblock {\em Asian J. Math.}, 3(1):93--107, 1999.
\newblock Sir Michael Atiyah: a great mathematician of the twentieth century.

\bibitem[MM22]{MM}
A.~Mori and Y.~Murakami.
\newblock Witten--{R}eshetikhin--{T}uraev invariants, homological blocks, and
  quantum modular forms for unimodular plumbing {H}-graphs.
\newblock {\em SIGMA. Symmetry, Integrability and Geometry: Methods and
  Applications}, 18:034, 2022.

\bibitem[MM24]{M-Mistegard}
W.~Misteg{\aa}rd and Y.~Murakami.
\newblock A proof of the radial limit conjecture for
  {C}ostantino--{G}eer--{P}atureau-{M}irand quantum invariants.
\newblock 2024.

\bibitem[MT21]{Matsusaka-Terashima}
T.~Matsusaka and Y.~Terashima.
\newblock Modular transformations of homological blocks for {S}eifert fibered
  homology $3$-spheres.
\newblock {\em arXiv:2112.06210}, 2021.

\bibitem[MT23]{M-Terashima}
Y.~Murakami and Y.~Terashima.
\newblock Homological blocks with simple {L}ie algebras and
  {W}itten--{R}eshetikhin--{T}uraev invariants.
\newblock {\em arXiv:2308.04010}, 2023.

\bibitem[Mur22a]{M:plumbed}
Y.~Murakami.
\newblock Witten--{R}eshetikhin--{T}uraev invariants and homological blocks for
  plumbed homology spheres.
\newblock 2022.
\newblock arXiv:2205.01282.

\bibitem[Mur22b]{M:indefinite}
Y.~Murakami.
\newblock Witten--{R}eshetikhin--{T}uraev invariants and indefinite false theta
  functions for plumbing indefinite {H}-graphs.
\newblock {\em to appear in Journal of Knot Theory and Its Ramifications,
  arXiv:2212.09972}, 2022.

\bibitem[Mur23]{M:GPPV}
Y.~Murakami.
\newblock A proof of a conjecture of {G}ukov--{P}ei--{P}utrov--{V}afa.
\newblock {\em to appear in Communications in Mathematical Physics,
  arXiv:2302.13526}, 2023.

\bibitem[Par20]{Park}
S.~Park.
\newblock Higher rank $\hat{Z}$ and $f_k$.
\newblock {\em SIGMA. Symmetry, Integrability and Geometry: Methods and
  Applications}, 16:044, 2020.

\bibitem[RT91]{Reshetikhin-Turaev}
N.~Reshetikhin and V.~G. Turaev.
\newblock Invariants of {$3$}-manifolds via link polynomials and quantum
  groups.
\newblock {\em Invent. Math.}, 103(3):547--597, 1991.

\bibitem[Wit89]{Witten}
E.~Witten.
\newblock Quantum field theory and the {J}ones polynomial.
\newblock {\em Comm. Math. Phys.}, 121(3):351--399, 1989.

\bibitem[Won01]{Wong}
R.~Wong.
\newblock {\em Asymptotic approximations of integrals}, volume~34 of {\em
  Classics in Applied Mathematics}.
\newblock Society for Industrial and Applied Mathematics (SIAM), Philadelphia,
  PA, 2001.
\newblock Corrected reprint of the 1989 original.

\bibitem[Zag06]{Zagier:asymptotic}
D.~Zagier.
\newblock The mellin transform and other useful analytic techniques.
\newblock In {\em Appendix to E. Zeidler, Quantum Field Theory {I}: {B}asics in
  Mathematics and Physics. {A} Bridge Between Mathematicians and Physicists},
  pages 305--323. Springer, Berlin, 2006.

\end{thebibliography}

% --------------------------------------------------------------------------
\end{document}